\documentclass[a4paper]{amsart}

\usepackage{listings}
\usepackage{pstricks, pst-node}

\newtheorem{thm}{Theorem}[section]

\newtheorem{conj}{Conjecture}[section]
\theoremstyle{definition}
\newtheorem*{maindef}{Definition}

\usepackage[colorlinks=true, linkcolor=blue, citecolor=blue,
  urlcolor=blue, pdfauthor={Dan Drake, Jang Soo Kim}]{hyperref}

\DeclareMathOperator{\dcr}{dcr}
\DeclareMathOperator{\Dcr}{DCR}
\DeclareMathOperator{\dne}{dne}
\DeclareMathOperator{\Dne}{DNE}

\DeclareMathOperator{\type}{type}

\DeclareMathOperator{\Ch}{Ch}

\def\sch.{Schr{\"o}der}

\def\O{\mathcal{O}}
\def\C{\mathcal{C}}
\def\S{\mathcal{S}} 
\def\T{\mathcal{T}}

\def\bk#1{\{#1\}}
\psset{unit=0.7cm, dotsize=5pt}
\psset{subgriddiv=1,gridcolor=lightgray,gridlabels=0pt,gridwidth=0.4pt}

\def\vput#1{\pnode(#1,1){#1} \pscircle*(#1,1){.1} \rput(#1,.5){$#1$}}
\def\edge#1#2{\ncarc[arcangle=50]{#1}{#2}}
\def\LOOP#1{\nccircle{#1}{0.15}}

\def\OPENER{\pnode(0,1){vertex} \pnode(0.25,1.3){half}
\ncarc[arcangle=-50]{half}{vertex} \pscircle*(0,1){.1}}

\def\opener#1{\rput(#1,0){\OPENER} \rput(#1,.5){$#1$}}

\newcount\sx \newcount\sy \newcount\ex \newcount\ey
\def\CHS#1(#2,#3){ 
\sx=#2 \sy=#3 \ex=#2 \ey=#3
\advance\ex by1 \advance\ey by0
\psline[linewidth=1.5pt](\sx,\sy)(\ex,\ey)
\rput(\number\sx.5,\number\sy.3){$#1$}
\psdot(\number\sx,\number\sy) \psdot(\number\ex,\number\ey)
}
\def\CDS#1(#2,#3){ 
\sx=#2 \sy=#3 \ex=#2 \ey=#3
\advance\ex by1 \advance\ey by-1
\psline[linewidth=1.5pt](\sx,\sy)(\ex,\ey)
\rput(\number\sx.7,\number\ey.7){$#1$}
\psdot(\number\sx,\number\sy) \psdot(\number\ex,\number\ey)
}
\def\HS(#1,#2){ 
\sx=#1 \sy=#2 \ex=#1 \ey=#2
\advance\ex by1 \advance\ey by0
\psline[linewidth=1.5pt](\sx,\sy)(\ex,\ey)
\psdot(\number\sx,\number\sy) \psdot(\number\ex,\number\ey)
}
\def\US(#1,#2){ 
\sx=#1 \sy=#2 \ex=#1 \ey=#2
\advance\ex by1 \advance\ey by1
\psline[linewidth=1.5pt](\sx,\sy)(\ex,\ey)
\psdot(\number\sx,\number\sy) \psdot(\number\ex,\number\ey)
}
\def\DS(#1,#2){ 
\sx=#1 \sy=#2 \ex=#1 \ey=#2
\advance\ex by1 \advance\ey by-1
\psline[linewidth=1.5pt](\sx,\sy)(\ex,\ey)
\psdot(\number\sx,\number\sy) \psdot(\number\ex,\number\ey)
}
\def\DHS(#1,#2){ 
\sx=#1 \sy=#2 \ex=#1 \ey=#2
\advance\ex by2 \advance\ey by0
\psline[linewidth=1.5pt](\sx,\sy)(\ex,\ey)
\psdot(\number\sx,\number\sy) \psdot(\number\ex,\number\ey)
}

\begin{document}

\title[$k$-distant crossings and nestings]{$k$-distant crossings and
  nestings of matchings and partitions}

\author{Dan Drake}
\address{Department of Mathematical Sciences\\
  Korea Advanced Institute of Science and Technology\\
  Daejeon, Korea} \email[Dan Drake]{ddrake@member.ams.org} \urladdr[Dan
Drake]{http://mathsci.kaist.ac.kr/~drake} \thanks{The first author was
  supported by the second stage of the Brain Korea 21 Project, The
  Development Project of Human Resources in Mathematics, KAIST in 2008.}

\author{Jang Soo Kim} \email[Jang Soo Kim]{jskim@kaist.ac.kr}
\urladdr[Jang Soo Kim]{http://combinat.kaist.ac.kr/~jskim}

\keywords{crossings, nestings, set partitions, matchings}
\subjclass[2000]{Primary: 05A15; Secondary: 05A18, 33C45, 05-04}

\date{\today}

\begin{abstract}
  We define and consider $k$-distant crossings and nestings for
  matchings and set partitions, which are a variation of crossings and
  nestings in which the distance between vertices is important. By
  modifying an involution of Kasraoui and Zeng (Electronic J.
  Combinatorics 2006, research paper 33), we show that the joint
  distribution of $k$-distant crossings and nestings is symmetric. We
  also study the numbers of $k$-distant noncrossing matchings and
  partitions for small $k$, which are counted by well-known sequences,
  as well as the orthogonal polynomials related to $k$-distant
  noncrossing matchings and partitions. We extend Chen et al.'s
  $r$-crossings and enhanced $r$-crossings.
\end{abstract}


\maketitle

\section{Introduction} \label{s:intro}

A \emph{(set) partition} of $[n]=\{1,2,\ldots,n\}$ is a set of disjoint
subsets of $[n]$ whose union is $[n]$. Each element of a partition is
called a \emph{block}. We will write a partition as a sequence of
blocks, for instance, $\bk{1,4,8} \bk{2,5,9} \bk{3} \bk{6,7}$. Let
$\Pi_n$ denote the set of partitions of $[n]$.

\begin{figure}[h]
  \begin{center}
    \begin{pspicture}(1,0.5)(12,2) \vput{1} \vput{2} \vput{3} \vput{4}
      \vput{5} \vput{6} \vput{7} \vput{8} \vput{9} \vput{10} \vput{11}
      \vput{12} \edge{1}{5} \edge{2}{4} \edge{4}{9} \LOOP{3}
      \edge{6}{12} \edge{7}{10} \edge{10}{11} \LOOP{8}
    \end{pspicture}
    \caption{Diagram for
      $\bk{1,5}\bk{2,4,9}\bk{3}\bk{6,12}\bk{7,10,11}\bk{8}$.}
    \label{f:genjiko-diagram}
  \end{center}
\end{figure}
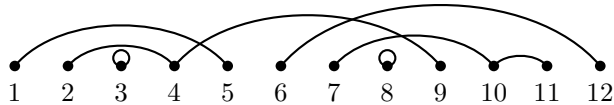

Let $\pi$ be a partition of $[n]$. A \emph{vertex} of $\pi$ is an
integer $i\in [n]$. An \emph{edge} of $\pi$ is a pair $(i,j)$ of
vertices satisfying either (1) $i<j$, and $i$ and $j$ are in the same
block with no vertex between them in that block, or (2) $i=j$ and the
block containing $i$ has no other vertex. Thus when we arrange vertices
of $\pi=\bk{1,5}\bk{2,4,9}\bk{3}\bk{6,12}\bk{7,10,11}\bk{8}$, in a line
in increasing order and draw edges we get \autoref{f:genjiko-diagram}.

A vertex $v$ of $\pi$ is called an \emph{opener} (resp. \emph{closer})
if $v$ is the smallest (resp. largest) element of a block consisting of
at least two integers. A vertex $v$ is called a \emph{singleton} if $v$
itself makes a block. A vertex $v$ is called a \emph{transient} if there
are two edges connected to $v$. Let $\O(\pi)$ (resp. $\C(\pi)$,
$\S(\pi)$, $\T(\pi)$) be the set of openers (resp. closers, singletons,
transients) of $\pi$. Let $\type(\pi)=(\O(\pi), \C(\pi), \S(\pi),
\T(\pi))$ and $\type'(\pi)=(\O(\pi), \C(\pi), \S(\pi) \cup \T(\pi))$.
For the partition in \autoref{f:genjiko-diagram}, the type of $\pi$ is
$\type(\pi) = (\bk{1,2,6,7}, \bk{5,9,11,12}, \bk{3,8}, \bk{4,10})$.

A \emph{(complete) matching} is a partition without singletons or
transients; this is the same thing as a partition in which all blocks
have size $2$.

Now we can define the main object of our study.

\begin{maindef}
  Let $k$ be a nonnegative integer. A \emph{$k$-distant crossing} of
  $\pi$ is a pair of edges $(i_1, j_1)$ and $(i_2, j_2)$ of $\pi$
  satisfying $i_1<i_2 \le j_1<j_2$ and $j_1 - i_2 \geq k$. A
  \emph{$k$-distant nesting} of $\pi$ is a set of two edges $(i_1, j_1)$
  and $(i_2, j_2)$ of $\pi$ satisfying $i_1<i_2\leq j_2<j_1$ and $j_2 -
  i_2 \geq k$.

  Let $\dcr_k(\pi)$ (resp. $\dne_k(\pi)$) denote the number of
  $k$-distant crossings ($k$-distant nestings) in $\pi$. Thus
  $\dcr_1(\pi)$ is the number of usual crossings of $\pi$.
\end{maindef}

For example, in the partition in \autoref{f:genjiko-diagram}, the edges
$(4,9)$ and $(6,12)$ form a $3$-distant crossing (as well as an
$i$-distant crossing for $i=0,1,2$), the edges $(1,5)$ and $(2,4)$ form
a $2$-distant nesting, the edges $(2,4)$ and $(4,9)$ form a $0$-distant
crossing, and the edges $(7,10)$ and $(8,8)$ form a $0$-distant nesting.
That partition has $\dcr_0(\pi) = 5$, $\dcr_2(\pi) = 2$, and
$\dne_2(\pi) = 2$.

Kasraoui and Zeng \cite{Kasraoui2006} found an involution $\varphi:\Pi_n
\rightarrow \Pi_n$ such that $\type(\varphi(\pi))=\type(\pi)$ and
$\dcr_1(\varphi(\pi))=\dne_1(\pi), \dne_1(\varphi(\pi))=\dcr_1(\pi)$.
Modifying this involution, for $k\geq0$, we find an involution
$\varphi_k : \Pi_n \rightarrow \Pi_n$ such that
$\dcr_k(\varphi_k(\pi))=\dne_k(\pi)$,
$\dne_k(\varphi_k(\pi))=\dcr_k(\pi)$ and
$\type(\varphi_k(\pi))=\type(\pi)$ if $k\geq1$;
$\type'(\varphi_k(\pi))=\type'(\pi)$ if $k=0$.

Noncrossing partitions and matchings are interesting and pervasive
objects that arise frequently in diverse areas of mathematics; see
\cite{McCammond2006} and \cite{Simion2000} and the references therein
for an introduction to noncrossing partitions. A partition $\pi$ is
called \emph{$k$-distant noncrossing} if $\pi$ has no $k$-distant
crossing. Let $NCM_k(n)$ denote the number of $k$-distant noncrossing
matchings of $[n]$. Let $NCP_k(n)$ denote the number of $k$-distant
noncrossing partitions of $[n]$.

\autoref{t:dnc-cm} and \autoref{t:dnc-setp} show $NCM_k(n)$ and
$NCP_k(n)$ for small values of $n$ and $k$. We use $k = \infty$ to
indicate that $i$-distant crossing is allowed for any positive integer
$i$, so that $NCM_\infty(n)$ and $NCP_\infty(n)$ equal the total number
of matchings of $[2n]$ and partitions of $[n]$, respectively. A matching
or partition cannot have a $k$-distant crossing for $k > n - 3$, so for
fixed $n$, $NCM_k(n)$ and $NCP_k(n)$ will ``converge'' to the number of
matchings and number of partitions, respectively; for readability we
omit those numbers in the tables. The $n=0$ column is all $1$'s for both
tables, of course.

\begin{table}
  \centering
  \begin{tabular}{r|rrrrrrrrrr}
    $k$ \textbackslash $n$
    & 2 & 4 &  6 &   8 &  10 &    12 &     14 &      16 &       18 &        20 \\
    \hline                                                  
    1 & 1 & 2 &  5 &  14 &  42 &   132 &    429 &    1430 &     4862 &     16796 \\
    2 &   & 3 & 11 &  45 & 197 &   903 &   4279 &   20793 &   103049 &    518859 \\
    3 &   &   & 14 &  71 & 387 &  2210 &  13053 &   79081 &   488728 &   3069007 \\
    4 &   &   & 15 &  91 & 581 &  3906 &  27189 &  194240 &  1416168 &  10494328 \\
    5 &   &   &    & 102 & 753 &  5752 &  45636 &  372360 &  3101523 &  26266917 \\
    6 &   &   &    & 105 & 873 &  7541 &  66690 &  607128 &  5657520 &  53631564 \\
    7 &   &   &    &     & 930 &  8985 &  88450 &  885394 &  9067611 &  94719138 \\
    8 &   &   &    &     & 945 &  9885 & 107847 & 1187376 & 13233511 & 150234570 \\
    9 &   &   &    &     &     & 10290 & 122115 & 1476948 & 17933348 & 219754737 \\
    10 &   &   &    &     &     & 10395 & 130515 & 1715475 & 22701570 & 300724081 \\
    11 &   &   &    &     &     &       & 134190 & 1881495 & 26969370 & 386669322 \\
    12 &   &   &    &     &     &       & 135135 & 1975995 & 30306045 & 468680940 \\
    13 &   &   &    &     &     &       &        & 2016630 & 32546745 & 538581120 \\
    14 &   &   &    &     &     &       &        & 2027025 & 33794145 & 591287445 \\
    15 &   &   &    &     &     &       &        &         & 34324290 & 625810185 \\
    16 &   &   &    &     &     &       &        &         & 34459425 & 652702050 \\
    17 &   &   &    &     &     &       &        &         &          & 644729085 \\
    18 &   &   &    &     &     &       &        &         &          & 654729075 \\
    \hline 
    $\infty$
    & 1 & 3 & 15 & 105 & 945 & 10395 & 135135 & 2027025 & 34459425 & 654729075 \\
  \end{tabular}                                              
  \caption{$k$-distant noncrossing matchings. The $k=0$ row is omitted
    because, as matchings have no transient vertices, the $k=0$ row is
    the same as $k=1$ row; both, of course, are counted by the Catalan
    numbers (A000108). The $k=2$ row is the little \sch. numbers
    (A001003).}
  \label{t:dnc-cm}
\end{table}

\begin{table}
  \centering
  \begin{tabular}{r|rrrrrrrrrrrr}
    $k$ \textbackslash $n$
    & 1 & 2 & 3 &  4 &  5 &   6 &   7 &    8 &     9 &     10 &     11 &      12\\
    \hline
    0 & 1 & 2 & 4 &  9 & 21 &  51 & 127 &  323 &   835 &   2188 &   5798 &   15511\\
    1 &   &   & 5 & 14 & 42 & 132 & 429 & 1430 &  4862 &  16796 &  58786 &  208012\\
    2 &   &   &   & 15 & 51 & 188 & 731 & 2950 & 12235 &  51822 & 223191 &  974427\\
    3 &   &   &   &    & 52 & 201 & 841 & 3726 & 17213 &  82047 & 400600 & 1993377\\
    4 &   &   &   &    &    & 203 & 872 & 4037 & 19796 & 101437 & 537691 & 2926663\\
    5 &   &   &   &    &    &     & 877 & 4125 & 20802 & 110950 & 618777 & 3575688\\
    6 &   &   &   &    &    &     &     & 4140 & 21095 & 114663 & 657698 & 3943294\\
    7 &   &   &   &    &    &     &     &      & 21147 & 115772 & 673019 & 4118232\\
    8 &   &   &   &    &    &     &     &      &       & 115975 & 677693 & 4187838\\
    9 &   &   &   &    &    &     &     &      &       &        & 678570 & 4209457\\
    10 &   &   &   &    &    &     &     &      &       &        &        & 4213597\\
    \hline
    $\infty$ & 1 & 2 & 5 & 15 & 52 & 203 & 877 & 4140 & 21147 & 115975 & 678570 & 4213597
  \end{tabular}
  \caption{$k$-distant noncrossing  partitions. The $k=0, 1$, and $2$ rows
    are counted by Motzkin numbers (A001006), the Catalan numbers, and
    A007317, respectively.}
  \label{t:dnc-setp}
\end{table}

It is well known that noncrossing matchings of $[2n]$ and noncrossing
partitions of $[n]$ are counted by the Catalan number $C_n$. Thus
$NCM_0(2n) = NCM_1(2n)=NCP_1(n)=C_n$. We will show that $NCM_2(2n) =
s_n$ and $NCP_0(n) = M_n$, where $s_n$ and $M_n$ are the little \sch.
numbers (A001003 in \cite{Sloane}) and the Motzkin numbers (A001006 in
\cite{Sloane}) respectively. We will also find the generating functions
for $NCP_2(n)$ and $NCM_3(2n)$.

Throughout this paper we will frequently refer to sequences in the
Online Encyclopedia of Integer Sequences \cite{Sloane} using their ``A
number''; we will usually omit the citation to \cite{Sloane} and
consider it understood that things like ``A000108'' are a reference to
the corresponding sequence in the OEIS.

The rest of this paper is organized as follows. In
\autoref{s:kas-zeng-inv}, we modify Kasraoui and Zeng's involution to
prove the joint distribution of $k$-distant crossings and nestings is
symmetric. In \autoref{s:motzkin-path-charlier-diagram}, we review a
bijection between partitions and Charlier diagrams. In
\autoref{s:k-dnc-matchings} and \autoref{s:k-dnc-partitions}, we study
the number of $k$-distant noncrossing matchings and partitions, and, in
\autoref{s:orthpoly}, we consider the orthogonal polynomials related to
these numbers. In \autoref{s:kr-dcr}, we extend $r$-crossings and
enhanced $r$-crossings of Chen et al. \cite{Chen2007}. We include an
appendix of Sage code used to compute the entries of \autoref{t:dnc-cm}
and \autoref{t:dnc-setp}.

\section{Modification of the involution of Kasraoui and Zeng}
\label{s:kas-zeng-inv}

Kasraoui and Zeng \cite{Kasraoui2006} found an involution $\varphi:\Pi_n
\rightarrow \Pi_n$ such that $\dcr_1(\varphi(\pi))=\dne_1(\pi)$,
$\dne_1(\varphi(\pi))=\dcr_1(\pi)$ and $\type(\varphi(\pi))=\type(\pi)$.
In this section, for fixed $k\geq0$, we find an involution
$\varphi_k:\Pi_n \rightarrow \Pi_n$ such that
$\dcr_k(\varphi_k(\pi))=\dne_k(\pi)$ and
$\dne_k(\varphi_k(\pi))=\dcr_k(\pi)$. Since complete matchings can be
thought of as set partitions with blocks all of size two, this
involution will also show that the distribution of $\dcr_{k}$ and
$\dne_{k}$ is symmetric.

We will follow Kasraoui and Zeng's notations. We will identify a
partition $\pi$ to its diagram as shown in \autoref{f:genjiko-diagram}.

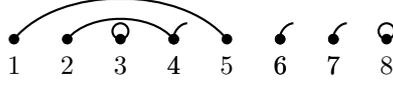
\begin{figure}
  \begin{center}
    \begin{pspicture}(1,0.5)(8,2) \vput{1} \vput{2} \vput{3} \vput{4}
      \vput{5} \vput{6} \vput{7} \vput{8}\edge{1}{5} \edge{2}{4}
      \opener{4} \LOOP{3} \opener{6} \opener{7} \LOOP{8}
    \end{pspicture}
  \end{center}
  \caption{The $8$-th trace $T_8(\pi)$ of $\pi$ in
    \autoref{f:genjiko-diagram}. The vacant vertices are $4,6$ and $7$.}
  \label{fig:trace}
\end{figure}

The \emph{$i$-th trace} $T_i(\pi)$ of $\pi$ is the diagram obtained from
$\pi$ by removing vertices greater than $i$. If a vertex $v\leq i$ is
connected to $u>i$ in $\pi$ then make a {\em half edge} from $v$ in
$T_i(\pi)$. Each vertex with a half edge is called {\em vacant} vertex.
For an example, see \autoref{fig:trace}.

Let $k$ be a fixed nonnegative integer. We define $\varphi_k:\Pi_n
\rightarrow \Pi_n$ as follows.

\begin{enumerate}
\item Set $T_0^{(k)} = \emptyset$.
\item For $1\leq i\leq n$, $T_i^{(k)}$ is obtained as follows.
  \begin{enumerate}
  \item Let $T_i^{(k)}$ (resp. $T_i'(\pi)$) be $T_{i-1}^{(k)}$ (resp.
    $T_{i-1}(\pi)$) with new vertex $i$.
  \item If $i\in\O(\pi)\cup\S(\pi)\cup\T(\pi)$, then make a half edge
    from $i$ both in $T_i^{(k)}$ and $T_i'(\pi)$.
  \item If $i\in\C(\pi)\cup\S(\pi)\cup\T(\pi)$, let $j$ be the vertex
    connected to $i$ in $\pi$.
    \begin{enumerate}
    \item If $i-j<k$, then $j$ must be a vacant vertex in $T_i^{(k)}$.
      Remove the half edge from $j$ and add an edge $(i,j)$ in
      $T_i^{(k)}$.
    \item If $i-j\geq k$, then let $U$ (resp. $V$) be the set of all
      vacant vertices $v$ in $T_i^{(k)}$ (resp. $T_i'(\pi)$) such that
      $i-v\geq k$. Let $\gamma_i^{(k)}(\pi)$ denote the integer $r$ such
      that $j$ is the $r$-th largest element of $V$. Let $j'$ be the
      $\gamma_i^{(k)}(\pi)$-th smallest element of $U$. Remove the half
      edge from $j'$ and add an edge $(j',i)$ in $T_i^{(k)}$.
    \end{enumerate}
  \end{enumerate}
\item Set $\varphi_k(\pi)=T_n^{(k)}$.
\end{enumerate}

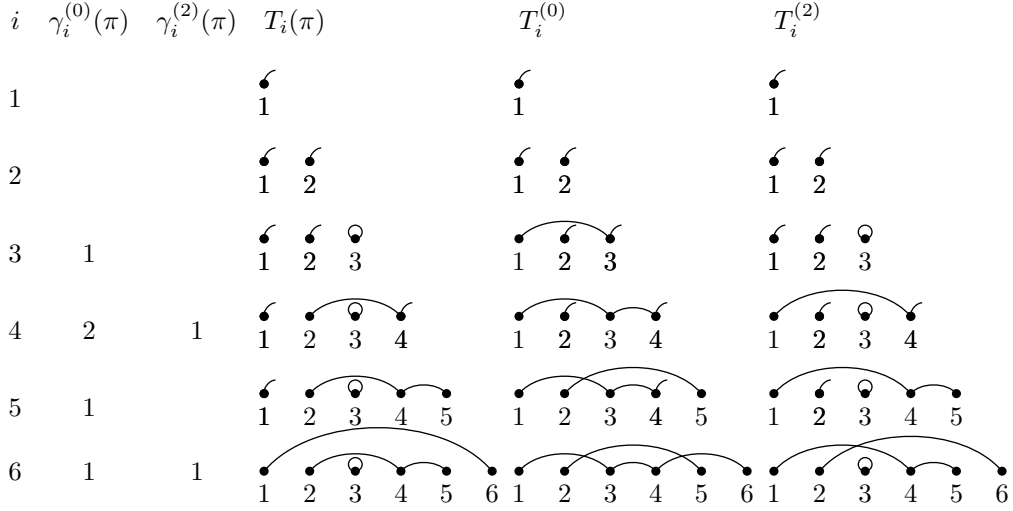
\begin{figure}
  \psset{unit=0.6cm, dotsize=3.5pt, linewidth=0.5pt}
  \begin{tabular}{ccclll}
    $i$ & $\gamma_i^{(0)}(\pi)$ & $\gamma_i^{(2)}(\pi)$ & $T_i(\pi)$ & $T^{(0)}_i$ & $T^{(2)}_i$ \\
    1 & & &
    \begin{pspicture}(1,0.5)(6,2)
      \vput{1} \opener{1}
    \end{pspicture}
    &    \begin{pspicture}(1,0.5)(6,2)
      \vput{1} \opener{1}
    \end{pspicture}
    &    \begin{pspicture}(1,0.5)(6,2)
      \vput{1} \opener{1}
    \end{pspicture}\\
    2 & & &
    \begin{pspicture}(1,0.5)(6,2)
      \vput{1} \opener{1} \vput{2} \opener{2}
    \end{pspicture}
    &    \begin{pspicture}(1,0.5)(6,2)
      \vput{1} \opener{1} \vput{2} \opener{2}
    \end{pspicture}
    &    \begin{pspicture}(1,0.5)(6,2)
      \vput{1} \opener{1} \vput{2} \opener{2}
    \end{pspicture}\\
    3 & 1& &
    \begin{pspicture}(1,0.5)(6,2)
      \vput{1} \opener{1} \vput{2} \opener{2} \vput3 \LOOP{3}
    \end{pspicture}
    &    \begin{pspicture}(1,0.5)(6,2)
      \vput{1} \vput{2} \vput{3} \opener{2} \opener{3} \edge{1}{3}
    \end{pspicture}
    &    \begin{pspicture}(1,0.5)(6,2)
      \vput{1} \opener{1} \vput{2} \opener{2} \vput3 \LOOP{3}
    \end{pspicture}\\
    4 & 2& 1&
    \begin{pspicture}(1,0.5)(6,2)
      \vput{1} \opener{1} \vput2 \vput3 \vput4 \edge24 \LOOP3 \opener4
    \end{pspicture}
    &    \begin{pspicture}(1,0.5)(6,2)
      \vput{1} \vput{2} \vput{3} \opener{2} \edge{1}{3} \vput4 \edge34 \opener4
    \end{pspicture}
    &    \begin{pspicture}(1,0.5)(6,2)
      \vput{1} \opener2 \vput2 \vput3 \vput4 \edge14 \LOOP3 \opener4
    \end{pspicture}\\
    5 & 1& &
    \begin{pspicture}(1,0.5)(6,2)
      \vput{1} \opener{1} \vput2 \vput3 \vput4 \edge24 \LOOP3 \vput5 \edge45 
    \end{pspicture}
    &    \begin{pspicture}(1,0.5)(6,2)
      \vput{1} \vput{2} \vput{3} \vput5 \edge{1}{3} \vput4 \edge34 \opener4 \edge25
    \end{pspicture}
    &    \begin{pspicture}(1,0.5)(6,2)
      \vput{1} \opener2 \vput2 \vput3 \vput4 \vput5 \edge14 \LOOP3 \edge45
    \end{pspicture}\\
    6 & 1& 1&
    \begin{pspicture}(1,0.8)(6,2)
      \vput{1} \vput2 \vput3 \vput4 \edge24 \LOOP3 \vput5 \edge45 \vput6 \edge16
    \end{pspicture}
    &    \begin{pspicture}(1,0.8)(6,2)
      \vput{1} \vput{2} \vput{3} \vput5 \vput6 \edge{1}{3} \vput4 \edge34 \edge46 \edge25
    \end{pspicture}
    &    \begin{pspicture}(1,0.8)(6,2)
      \vput{1} \vput2 \vput3 \vput4 \vput5 \vput6 \edge14 \LOOP3 \edge45 \edge26
    \end{pspicture}\\
  \end{tabular}
  \caption{Construction of $\varphi_0(\pi)=T_6^{(0)}$ and
    $\varphi_2(\pi) = T_6^{(2)}$ for $\pi=\bk{1,6}\bk{2,4,5}\bk{3}$.}
  \label{fig:construction}
\end{figure}

For example, see \autoref{fig:construction}. Using the same argument as
in \cite{Kasraoui2006}, we can prove that $\varphi_k$ is an involution
and satisfies $\dcr_k(\varphi_k(\pi))=\dne_k(\pi)$,
$\dne_k(\varphi_k(\pi))=\dcr_k(\pi)$, $\type(\varphi_k(\pi))=\type(\pi)$
if $k\geq1$; $\type'(\varphi_k(\pi))=\type'(\pi)$ if $k=0$. Thus we have
the following.

\begin{thm}
  Let $k$ be a nonnegative integer. Then
  \[ \sum_{\pi\in\Pi_n} x^{\dcr_k(\pi)} y^{\dne_k(\pi)} =
  \sum_{\pi\in\Pi_n} x^{\dne_k(\pi)} y^{\dcr_k(\pi)}. \]
\end{thm}

\section{Motzkin paths and Charlier diagrams}
\label{s:motzkin-path-charlier-diagram}

In this section, we recall a bijection between partitions and Charlier
diagrams \cite{Flajolet1980, Kasraoui2006}.

A \emph{step} is a pair $(p,q)$ of points with $p$ and $q$ in $\mathbb Z
\times \mathbb Z$. The \emph{height} of a step $(p,q)$ is the second
component of $p$, i.e, if $p=(a,b)$ then the height of the step $(p,q)$
is~$b$. A step $(p,q)$ is called an \emph{up} (resp. \emph{down},
\emph{horizontal}) step if the component-wise difference $q-p$ is
$(1,1)$ (resp. $(1,-1)$, $(1,0)$). A \emph{path} of length $n$ is a
sequence $(p_0,p_1,p_2,\ldots,p_n)$ of $n+1$ points in $\mathbb Z \times
\mathbb Z$. The \emph{$i$-th step} of a path $(p_0,p_1,p_2,\ldots,p_n)$
is $(p_{i-1},p_i)$. A \emph{Motzkin path} of length $n$ is a path from
$(0,0)$ to $(n,0)$ consisting of up steps, down steps, and horizontal
steps that never goes below the $x$-axis. A \emph{Charlier diagram} of
length $n$ is a pair $(M, e)$ where $M=(p_0,p_1,\ldots,p_n)$ is a
Motzkin path of length $n$ and $e=(e_1,e_2,\ldots,e_n)$ is a sequence of
integers such that:
\begin{enumerate}
\item if the $i$-th step is an up step then $e_i=0$,
\item if the $i$-th step is a down step of height $h$ then $1\leq
  e_i\leq h$,
\item if the $i$-th step is a horizontal step of height $h$ then $0\leq
  e_i\leq h$.
\end{enumerate}

We will identify a Charlier diagram $(M,e)$ with the sequence
$(s_1,s_2,\ldots,s_n)$ of labeled letters in
$\{U,D_1,D_2,\ldots,H_0,H_1,H_2,\ldots\}$ such that $s_i=U$ (resp.
$s_i=D_{e_i}$, $s_i=H_{e_i}$) if the $i$-th step of $M$ is an up (resp.
down, horizontal) step.

Let $\pi$ be a partition of $[n]$. Recall that in the previous section,
if $i$ is a closer or transient, then $\gamma_i^{(1)}(\pi)$ is the
integer $r$ such that $i$ is connected to the $r$-th largest integer in
$T_{i-1}^{(1)}(\pi)$.

The corresponding Charlier diagram $\Ch(\pi)=(s_1,s_2,\ldots,s_n)$ is
defined as follows:
\begin{enumerate}
\item if $i$ is an opener in $\pi$ then $s_i=U$,
\item if $i$ is a closer in $\pi$ and $\gamma_i^{(1)}(\pi)=r$ then
  $s_i=D_r$,
\item if $i$ is a singleton in $\pi$ then $s_i=H_0$,
\item if $i$ is a transient in $\pi$ and $\gamma_i^{(1)}(\pi)=r$ then
  $s_i=H_r$.
\end{enumerate}

\begin{figure}
  \begin{center}
    \begin{pspicture}(0,0)(12,4) \psgrid(0,0)(12,3) \US(0,0) \US(1,1)
      \CHS{0}(2,2) \CHS{1}(3,2) \CDS{2}(4,2) \US(5,1) \US(6,2)
      \CHS{0}(7,3) \CDS{3}(8,3) \CHS{1}(9,2) \CDS{1}(10,2) \CDS{1}(11,1)
    \end{pspicture}
    \caption{The Charlier diagram for the partition of
      \autoref{f:genjiko-diagram}. The label $e_i$ is written above the
      horizontal and down steps.}
    \label{f:charlier-diagram}
  \end{center}
\end{figure}
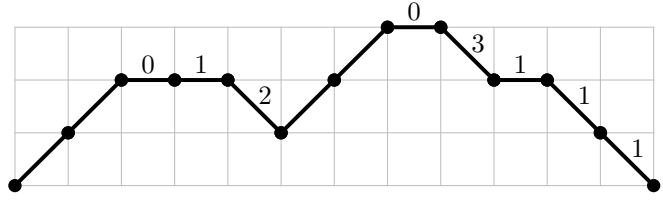

For example, see \autoref{f:charlier-diagram}.

It is easy to see that if there is a step $D_{\ell}$ or $H_{\ell}$ with
$\ell\geq 2$ in $\Ch(\pi)$, than $\pi$ has an $(\ell-1)$-distant
crossing.

\section{$k$-distant noncrossing matchings} \label{s:k-dnc-matchings}

In this section we will find the number of $k$-distant noncrossing
matchings for $k=0,1,2$ and $3$. Note that since there is no matching of
$[2n+1]$ we have $NCM_k(2n+1)=0$ for all $n$ and $k$. Thus we will only
consider $NCM_k(2n)$.

\subsection{$0$- and $1$-distant noncrossing matchings}
Since matchings have no transient vertices, being $0$-distant crossing
is equivalent to being $1$-distant crossing.

We can easily see that a matching $\pi$ is $1$-distant noncrossing if
and only if $\Ch(\pi)$ consists of $U$ and $D_1$. Thus a $1$-distant
noncrossing matching corresponds to a Dyck path.

\begin{thm}
  We have
  \[
  NCM_0(2n) = NCM_1(2n) =C_n = \frac{1}{n+1}\binom{2n}{n}.
  \]
\end{thm}

\subsection{$2$-distant noncrossing matchings}
Let $\pi$ be a $2$-distant noncrossing matching. Then $\Ch(\pi)$
consists of $U$, $D_1$, and $D_2$. By definition of $\Ch(\pi)$, $D_2$ is
of height at least $2$. Moreover, since $\pi$ has no $2$-distant
crossing, $D_2$ must immediately follow $U$. Thus we can consider
$\Ch(\pi)$ as a nonnegative path consisting of the three steps
$U=(1,1)$, $D_1=(1,-1)$ and $UD_2=(2,0)$ such that $UD_2$ never touches
the $x$-axis. This is exactly the definition of a \emph{little \sch.
  path}, see \cite{EC2}. Thus we get the following theorem.

\begin{thm}
  We have
  \[
  NCM_2(2n) = s_n,
  \]
  where $s_n$ is the little \sch. number (A001003).
\end{thm}

\subsection{$3$-distant noncrossing matchings}
\label{sec:3-dist-nc-matchings}

Let $\pi$ be a $3$-distant noncrossing matching. One can check that
$\Ch(\pi)$ consists of $U$, $D_1$, $D_2$, and $D_3$ satisfying the
following.
\begin{enumerate}
\item $D_\ell$ is of height at least $\ell$ for $\ell=1,2,3$.
\item $D_3$ can only occur after two consecutive $U$, and
\item $D_2$ can only occur after $U$ or after either $D_2$ or $D_3$
  which follows $U$.
\end{enumerate}

Thus we can consider $\Ch(\pi)$ as a path consisting of $6$ kinds of
steps: $U$ and $D_1$, which can appear at any height; and $UD_2$,
$UUD_3$, $UD_2D_2$, and $UUD_3D_2$, which can only appear above the line
$y=1$. Let $g(n)$ be the number of labeled Motzkin paths of length $n$
consisting of $U$, $D_1$, $UD_2$, $UUD_3$, $UD_2D_2$, and $UUD_3D_2$
with no restriction on height---these are not properly Charlier diagrams
since, for example, down steps at height $1$ may have label $2$. Let
$F(x)=\sum_{n\geq 0} NCM_3(2n)x^n$ and $G(x)=\sum_{n\geq 0} g(2n)x^n$.

Decomposing Charlier diagrams as in \autoref{tab:3dncmatch-path-decomp},
we can easily see that $G(x) = 1+(x+x^2)G(x)+(x^{1/2}+x^{3/2})^2G(x)^2$,
which means
\[
G(x)=\frac{1-x-x^2-\sqrt{1-6x-9x^2-2x^3+x^4}}{2x(x+1)^2}.
\]
Since $F(x) = 1 + x G(x) F(x)$, we get the generating function for
$NCM_3(2n)$.
\begin{thm}
  We have
  \begin{align*}
    \sum_{n\geq0} NCM_3(2n)x^n &=
    \frac{2(x+1)^2}{1+5x+3x^2+\sqrt{1-6x-9x^2-2x^3+x^4}}\\
    &=1+x+3\,{x}^{2}+14\,{x}^{3}+71\,{x}^{4}+387\,{x}^{5}+2210\,{x}^{6}+
    13053\,{x}^{7} +\cdots.
  \end{align*}
\end{thm}

\begin{table}
  \centering
  \begin{tabular}{l|l}
    Path type & Weight \\
    \hline
    $UD_{2}$ followed by any path & $x G(x)$ \\
    $UUD_{3}D_{2}$ followed by any path & $x^2 G(x)$ \\
    $U,$ any path, $D_{1}$, any path & $x G(x)^{2}$ \\
    $UUD_{3}$, any path, $D_{1}$, any path & $x^{2} G(x)^{2}$ \\
    $U,$ any path, $UD_{2}D_{2}$, any path & $x^{2}G(x)^{2}$ \\
    $UUD_{3}$, any path, $UD_{2}D_{2}$, any path & $x^{3} G(x)^{2}$
  \end{tabular}
  \caption{The six possible first-return decompositions for non-empty
    paths counted by the generating function $G(x)$ of
    \autoref{sec:3-dist-nc-matchings}.}
  \label{tab:3dncmatch-path-decomp}
\end{table}

\section{$k$-distant noncrossing partitions} \label{s:k-dnc-partitions}

\subsection{$0$-distant noncrossing partitions}
Let $\pi$ be a $0$-distant noncrossing partition. Then $\Ch(\pi)$
consists of $U, D_1, H_0$. Thus $\Ch(\pi)$ is a Motzkin path.
\begin{thm}
  The number of $0$-distant noncrossing partitions of $[n]$ is equal to
  the number of Motzkin paths of length $[n]$ (A001006).
\end{thm}

\subsection{$1$-distant noncrossing partitions}
Let $\pi$ be a $1$-distant noncrossing partition. Then $\pi$ is a usual
noncrossing partition. It is well known that the number of noncrossing
partitions of $[n]$ is the Catalan number $C_n$.

\begin{thm}
  We have
  \[
  NCP_1(n) = C_n.
  \]
\end{thm}

\subsection{$2$-distant noncrossing partitions}

Let $\pi$ be a $2$-distant noncrossing partition. Then $\Ch(\pi)$
consists of $U$, $D_1$, $D_2$, $H_0$, $H_1$, and $H_2$ and satisfies
\begin{enumerate}
\item $D_{\ell}$ and $H_{\ell}$ are of height at least $\ell$,
\item $H_2$ and $D_2$ can only occur after $U$, $H_1$, or $H_2$.
\end{enumerate}

Thus we can consider $\Ch(\pi)$ as a path with the following steps:
$UH_2^k$, $UH_2^kD_2$, $H_1H_2^k$, $H_1H_2^kD_2$, $H_0$, and $D_1$,
where $k$ is a nonnegative integer and $H_2^k$ means $k$ consecutive
$H_{2}$ steps.

Let $a(n)$ (resp. $b(n)$) denote the number of Charlier diagrams of
length $n$ consisting of the above steps such that $D_\ell$ and $H_\ell$
is of height at least $\ell-2$ (resp. at least $\ell-1$). In fact, the
height condition is unnecessary for $a(n)$ since every step is of height
at least $0$. Let $F(x)=\sum_{n\geq0} NCP_2(n) x^n$, $A(x)=\sum_{n\geq0}
a(n) x^n$, and $B(x)=\sum_{n\geq0} b(n) x^n$.

Note that the steps which increase the $y$-coordinate by $1$ are
$UH_2^k$; the steps which do not change the $y$-coordinate are $H_0$,
$H_1H_2^k$, and $UH_2^kD_2$; and the steps which decrease $y$-coordinate
by $1$ are $D_1$ and $H_1H_2^kD_2$. By decomposing Charlier diagrams as
we did with $G(x)$ in \autoref{sec:3-dist-nc-matchings}, we get
\begin{align*}
  A(x) &= 1+\left(x+\frac x{1-x}+\frac{x^2}{1-x}\right)A(x)
  +\frac{x}{1-x}
  \cdot \left(x+\frac{x^2}{1-x}\right) A(x)^2,\\
  B(x) &= 1+\left(2x+\frac{x^2}{1-x}\right)B(x) +\frac{x}{1-x}
  \cdot \left(x+\frac{x^2}{1-x}\right) A(x)B(x), \text{ and}\\
  F(x) &= 1+xF(x) + x^2 B(x)F(x).
\end{align*}
Solving these equations, we get the following theorem.
\begin{thm}\label{thm:ncp2}
  We have
  \begin{align*}
    \sum_{n\geq0} NCP_2(n) x^n &= \frac{3-3x-\sqrt{1-6x+5x^2}}{2(1-x)}=
    \frac 32 -\frac12 \sqrt{\frac{1-5x}{1-x}}\\
    &=1+x+2\,{x}^{2}+5\,{x}^{3}+15\,{x}^{4}+51\,{x}^{5}+188\,{x}^{6}+731\,{x
    }^{7}+2950\,{x}^{8}+\cdots.
  \end{align*}
\end{thm}

This sequence is A007317. Mansour and Severini \cite{Mansour2007} proved
that the generating function for the number of $12312$-avoiding
partitions is equal to that in \autoref{thm:ncp2}. Thus the number of
$2$-distant noncrossing partitions of $[n]$ is equal to the number of
$12312$-avoiding partitions of $[n]$. Yan \cite{Yan2008} found a
bijection from $12312$-avoiding partitions of $[n]$ to UH-free \sch.
paths of length $2n-2$. Composing several bijections including Yan's
bijection, Kim \cite{Kim2008pre} found a bijection between $2$-distant
noncrossing partitions and $12312$-avoiding partitions.

\section{Orthogonal polynomials} \label{s:orthpoly}

Given a sequence $\bk{\mu_n}_{n \ge 0}$, one may try to define a
sequence of polynomials $\bk{P_n(x)}_{n \ge 0}$ that are orthogonal with
respect to $\bk{\mu_n}$; that is, if we define a measure with $\mu_n =
\int x^n \,\mathrm{d}\mu$, then
\[
\int P_n(x) P_m(x) \,\mathrm{d}\mu = 0
\]
whenever $n \ne m$. Any sequence of polynomials satisfying the above
orthogonality relation must satisfy a three-term recurrence relation of
the form
\begin{equation}
  P_{n+1}(x) = (x - b_n) P_n(x) - \lambda_n P_{n-1}(x),
  \label{e:op-recursion}
\end{equation}
with $P_0(x) = 1$ and $P_1(x) = x - b_0$. Viennot showed
\cite{Viennot1983, Viennot1985} that for any sequence
$\bk{\mu_n}$---which are called the \emph{moments}---one can interpret
the moment $\mu_n$ as the generating function for weighted Motzkin paths
of length $n$ in which up steps have weight $1$, horizontal steps of
height $k$ have weight $b_k$, and down steps of height $k$ have weight
$\lambda_k$; then the polynomials in \eqref{e:op-recursion} will be
orthogonal with respect to $\bk{\mu_n}_{n \ge 0}$.

Many classical combinatorial sequences have been interpreted as the
moment sequences for a set of orthogonal polynomials, and the
corresponding orthogonality relation proved with a sign-reversing
involution. In particular, it is known that:

\begin{itemize}
\item If $\mu_{2n+1} = 0$ and $\mu_{2n} = C_n$, the Catalan number, then
  $b_n = 0$ and $\lambda_n = 1$; the corresponding polynomials are
  Chebyshev polynomials of the second kind \cite{deSainteCatherine1984},
  which may be defined by $U_{n+1}(x) = x U_n (x) - U_{n-1}(x)$, with
  $U_0(x) = 1$ and $U_1(x) = x$. These moments are $NCM_0(n)$ (and
  $NCM_{1}(n)$).
\item If $\mu_{2n+1} = 0$ and $\mu_{2n} = (2n-1)!!$, then $b_n = 0$ and
  $\lambda_n = n$; the corresponding polynomials are Hermite polynomials
  \cite{Viennot1983}. These moments are $NCM_{\infty}(n)$.
\item If $\mu_n = M_n$, the $n$-th Motzkin number, then $b_n = 1$,
  $\lambda_n = 1$; the corresponding polynomials are shifted Chebyshev
  polynomials of the second kind: $U_n(x-1)$. See \cite[section~
  4.1]{Drake2006}. These moments are $NCP_0(n)$.
\item If $\mu_n = B_n$, the number of partitions of $[n]$, then $b_n =
  n+1$ and $\lambda_n = n$; the corresponding polynomials are Charlier
  polynomials (with $a = 1$) \cite{Viennot1983}. These moments are
  $NCP_{\infty}(n)$.
\end{itemize}

With these observations in mind, it is natural to try to use, say,
$NCM_k(n)$ as a sequence of moments. Letting $k$ go from $0$ to infinity
would then allow us to interpolate between Chebyshev polynomials and
Hermite polynomials; using $NCP_k(n)$ would give the corresponding
interpolation between shifted Chebyshev and Charlier polynomials.

What happens if we use $NCM_2(n)$ for the moments? We know that
$NCM_{2}(2n+1)$ is zero, and $NCM_{2}(2n)$ equals the little \sch.
number $s_{n}$, which means the corresponding sequence of $b_{n}$'s is
all zeros. We need only find the $\lambda_{n}$'s.

\begin{thm}
  \label{thm:ncm2-lambdas}
  If $\lambda_{2n+1} = 1$ and $\lambda_{2n} = 2$, then the corresponding
  weighted Dyck paths are counted by the little \sch. numbers.
\end{thm}

\begin{proof}
  Both sequences have the same generating function: if we weight the
  upsteps and downsteps of little \sch. paths by $x$ and double
  horizontal steps by $x^{2}$, then the generating function of such
  paths is
  \[
  \frac{1 + x^{2} - \sqrt{x^{4} - 6x^{2} + 1}}{4x^{2}}.
  \]
  See Stanley \cite[p. 178]{EC2}. On the other hand, consider the
  following weightings for Dyck paths:
  \[
  \lambda_{n} =
  \begin{cases}
    1 & \text{$n$ odd}\\
    2 & \text{$n$ even},
  \end{cases}
  \quad\text{and}\quad \lambda_{n} =
  \begin{cases}
    2 & \text{$n$ odd}\\
    1 & \text{$n$ even},
  \end{cases}
  \]
  with upsteps all weighted $1$. Let $A(x)$ and $B(x)$ respectively
  denote the generating functions of such paths. By decomposing paths by
  their first return to the $x$-axis, we have
  \[
  A(x) = \frac{1}{1 - x^{2} B(x)} \quad\text{and}\quad B(x) = \frac{1}{1
    - 2x^{2} A(x)};
  \]
  by substituting the expression for $B(x)$ into that for $A(x)$ and
  solving, we find that
  \[
  A(x) = \frac{1 + x^{2} \pm \sqrt{x^{4} - 6x^{2} + 1}}{4x^{2}}.
  \]
  Using ``$-$'' yields the correct generating function, which coincides
  with the known generating function for the little \sch. numbers.
\end{proof}

This means that the orthogonal polynomials corresponding to $NCM_{2}(n$)
have $b_n = 0$, $\lambda_{2n+1} = 1$, and $\lambda_{2n} = 2$. They are a
special case of polynomials studied by Kim and Zeng \cite{Kim2003}: use
$U_n(x, 2)$ in their paper. Vauchassade de Chaumont and Viennot
\cite{Vauchassade1984} also studied these polynomials, although they use
a different normalization for the moments and instead get the big \sch.
numbers.

If we attempt to do the same with $NCM_3(n)$, we get stuck: since
$NCM_3(2n+1) = 0$, we know that $b_n = 0$, but the $\lambda_n$ sequence
starts with
\begin{equation}
  1, 2, \frac{5}{2}, \frac{3}{10}, \frac{76}{5},
  -\frac{680}{57}, -\frac{2311}{7752}, \frac{1246001}{314296},
  \frac{114710016}{151553069},\dots .
  \label{e:badlambdas}
\end{equation}
Not only are some $\lambda_n$'s fractions, but some are negative, which
means prospects for polynomials with nice combinatorics are dim.

Let us try the same line of attack with $k$-distant noncrossing
partitions. Using $NCP_1(n)$---Catalan numbers---for a set of moments,
we get a shifted version of Chebyshev polynomials of the second kind:
$b_0 = 1$, all other $b_n = 2$, and all $\lambda_n = 1$. These
polynomials can be written $U_n(x-2)$, with slightly different initial
conditions: $U_0(x) = 1$ and $U_1(x) = x-1$. The easiest way to see why
these recurrence coefficients and initial conditions are orthogonal with
respect to the Catalan numbers is with a bijection between Motzkin paths
of length $n$ with the above weighting and Dyck paths of length $2n$:
take each up step $U$ and make it $UU$, take each down step $D$ and make
it $DD$, and take each horizontal step $H$ and make it either $UD$ or
$DU$---except for the horizontal step at height zero, which can only be
made into $UD$. This process turns a weighted Motzkin path of length $n$
into a Dyck path of length $2n$ and is easily shown to be a bijection.

When using $NCP_2(n)$ and $NCP_3(n)$ as the moments, we again get some
fractional coefficients, but they seem much nicer. We have computed the
following with Maple: if $\mu_n=NCP_2(n)$ then
\begin{align*}
  \{b_n\}_{n \geq 0} &= \left\{1,3-1,3-\frac{1}{2}, 3-\frac{1}{10},
    3-\frac{1}{65}, 3-\frac{1}{442}, 3-\frac{1}{3026}, \ldots\right\}
  \text{and}\\
  \{\lambda_n\}_{n\geq1} &= \left\{ 1, 1+1, 1+\frac{1}{4},
    1+\frac{1}{25},1+\frac{1}{169}, 1+\frac{1}{1156},
    1+\frac{1}{7921},\ldots \right\};
\end{align*}
if $\mu_n=NCP_3(n)$ then
\[
\{b_n\}_{n \geq0} = \left\{ 1,2,3,3,3, \ldots\right\} \quad \text{and}
\quad \{\lambda_n\}_{n \geq1} = \left\{1,2,2,2,2,\ldots\right\}.
\]
The first case is very interesting. The sequences of denominators of
$b_n$'s and $\lambda_n$'s appear in A064170 and A081068 respectively.
Based on the above evidence, we make the following conjecture.

\begin{conj}
  If $\mu_n=NCP_2(n)$ then $b_0=b_1=\lambda_1=1$, and for $n\geq2$
  \[
  b_n = 3-\frac{1}{F_{2n-1}F_{2n-3}} \quad \text{and} \quad \lambda_n =
  1+\frac1{(F_{2n-3})^2},
  \]
  where $F_n$ is the $n$-th Fibonacci number, i.e.,
  $F_{n+1}=F_n+F_{n-1}$ and $F_1=F_2=1$.

  If $\mu_n=NCP_3(n)$ then $b_0 = b_1 = \lambda_1 = 1, \lambda_2 = 2$,
  and, for $n \ge 3$, $b_n=3$ and $\lambda_n=2$.
\end{conj}

\section{$k$-distant $r$-crossing} \label{s:kr-dcr}

Chen et al. \cite{Chen2007} considered a different kind of crossing
number. Our definition of $k$-distant crossing can be applied to their
definition.

Let $k\geq 0$ and $r\geq 2$ be integers. A \emph{$k$-distant
  $r$-crossing} is a set of $r$ edges
$(i_1,j_1),(i_2,j_2),\ldots,(i_r,j_r)$ such that $i_1<i_2<\cdots<i_r\leq
j_1<j_2<\cdots<j_r$ and $j_1-i_r \geq k$. Similarly, a \emph{$k$-distant
  $r$-nesting} is a set of $r$ edges
$(i_1,j_1),(i_2,j_2),\ldots,(i_r,j_r)$ such that $i_1<i_2<\cdots<i_r\leq
j_r<j_{r-1}<\cdots<j_1$ and $j_r-i_r \geq k$. In \cite{Chen2007}, they
defined an $r$-crossing and an enhanced $r$-crossing, which are a
$1$-distant $r$-crossing and a $0$-distant $r$-crossing respectively.

Let $\Dcr_k(\pi)$ (resp. $\Dne_k(\pi)$) be the maximal $r$ such that
$\pi$ has a $k$-distant $r$-crossing (resp. $k$-distant $r$-nesting).
Let $f_{n,S,T}(k;i,j)$ denote the number of partitions $\pi$ of $[n]$
such that $\Dcr_k(\pi)=i$, $\Dne_k(\pi)=j$, $\O(\pi)=S$ and $\C(\pi)=T$.
Chen et al. \cite{Chen2007} proved that
$f_{n,S,T}(k;i,j)=f_{n,S,T}(k;j,i)$ for $k=0,1$. Krattenthaler
\cite{Krattenthaler2006} extended this result using growth diagrams.

Using Krattenthaler's growth diagram method, we can get the following
theorem.
\begin{thm}
  Let $n\geq1$ and $k\geq0$ be integers. Then
  \[
  f_{n,S,T}(k;i,j) = f_{n,S,T}(k;j,i).
  \]
\end{thm}

\section*{Appendix: Sage code}

In this appendix we provide code used to compute the values in
\autoref{t:dnc-cm} and \autoref{t:dnc-setp}. This code is for use with
the free open-source computer mathematics system Sage
(\url{http://sagemath.org}). The source code is available as a separate
file: select ``source'' from the ``other formats'' link on the abstract
page for this preprint
(\href{http://arxiv.org/abs/0812.2725}{arxiv.org/abs/0812.2725}), and
extract \texttt{sage-code-appendix.sage} from the file you download.

\lstset{language=Python, showstringspaces=false,
  basicstyle=\footnotesize\ttfamily, columns=flexible}

\begin{lstlisting}
r"""
Sage code for computing k-distant crossing numbers.

This code accompanies the article arxiv:0812.2725; see
http://arxiv.org/abs/0812.2725.

Right now, this code only computes k-dcrossings. If you are only
interested in the distribution, this is good enough because the extended
Kasraoui-Zeng involution tells us the distribution of k-dcrossings and
k-dnestings is symmetric. It would be nice, though, to have a function
which actually performed that involution.

AUTHORS: 
    -- Dan Drake (2008-12-15): initial version.

EXAMPLES:

The example given in the paper. Note that in this format, we omit fixed
points since they cannot create any sort of crossing.

    sage: dcrossing([(1,5), (2,4), (4,9), (6,12), (7,10), (10,11)])
    3

"""

#*****************************************************************************
# Copyright (C) 2008 Dan Drake <ddrake@member.ams.org>
#
# This program is free software: you can redistribute it and/or modify
# it under the terms of the GNU General Public License as published by
# the Free Software Foundation, either version 2 of the License, or (at
# your option) any later version.
#
# See http://www.gnu.org/licenses/.
#*****************************************************************************

def CompleteMatchings(n):
    """Return a generator for the complete matchings of the set [1..n].
    
    INPUT:
        n -- nonnegative integer

    OUTPUT:
        A generator for the complete matchings of the set [1..n], or,
        what is basically the same thing, complete matchings of the
        graph K_n. Each complete matching is represented by a list of
        2-element tuples.

    EXAMPLES:
    There are 3 complete matchings on 4 vertices:
        sage: [m for m in CompleteMatchings(4)]
        [[(3, 4), (1, 2)], [(2, 4), (1, 3)], [(2, 3), (1, 4)]]

    There are no complete matchings on an odd number of vertices; the
    number of complete matchings on an even number of vertices is a
    double factorial:
        sage: [len([m for m in CompleteMatchings(n)]) for n in [0..8]]
        [1, 0, 1, 0, 3, 0, 15, 0, 105]

    The exact behavior of CompleteMatchings(n) if n is not a nonnegative
    integer depends on what [1..n] returns, and also on what range(1,
    len([1..n])) is.

    """
    for m in matchingsset([1..n]): yield m
 
def matchingsset(L):
    """Return a generator for complete matchings of the sequence L.

    This is not really meant to be called directly, but rather by
    CompleteMatchings().

    INPUT:
        L -- a sequence. Lists, tuples, et cetera; anything that
        supports len() and slicing should work.

    OUTPUT:
        A generator for complete matchings on K_n, where n is the length
        of L and vertices are labeled by elements of L. Each matching is
        represented by a list of 2-element tuples.

    EXAMPLES:
        sage: [m for m in matchingsset(('a', 'b', 'c', 'd'))]
        [[('c', 'd'), ('a', 'b')], [('b', 'd'), ('a', 'c')], [('b', 'c'), ('a', 'd')]]

        There's only one matching of the empty set/list/tuple: the empty
        matching.

        sage: [m for m in matchingsset(())]
        [[]]
    """
    if len(L) == 0:
        yield []
    else:
        for k in range(1, len(L)):
            for m in matchingsset(L[1:k] + L[k+1:]):
                yield m + [(L[0], L[k])]

def dcrossing(m_):
    """Return the largest k for which the given matching or set
    partition has a k-distant crossing.

    INPUT:
       m -- a matching or set partition, as a list of 2-element tuples
       representing the edges. You'll need to call setp_to_edges() on
       the objects returned by SetPartitions() to put them into the
       proper format.

    OUTPUT:
       The largest k for which the object has a k-distant crossing.
       Matchings and set partitions with no crossings at all yield -1.

    EXAMPLES:
    The main example from the paper:
        sage: dcrossing(setp_to_edges(Set(map(Set, [[1,5],[2,4,9],[3],[6,12],[7,10,11],[8]]))))
        3

    A matching example:

        sage: dcrossing([(4, 7), (3, 6), (2, 5), (1, 8)])
        2

    TESTS:
    The empty matching and set partition are noncrossing:
        sage: dcrossing([])
        -1
        sage: dcrossing(Set([]))
        -1

    One edge:
        sage: dcrossing([Set((1,2))])
        -1
        sage: dcrossing(Set([Set((1,2))]))
        -1

    Set partition with block of size >= 3 is always at least
    0-dcrossing:
        sage: dcrossing(setp_to_edges(Set([Set((1,2,3))])))
        0
    """
    d = -1
    m = list(m_)
    while len(m) > 0:
        e1_ = m.pop()
        for e2_ in m:
            e1, e2 = sorted(e1_), sorted(e2_)
            if (e1[0] < e2[0] and e2[0] <= e1[1] and e1[1] < e2[1] and
                e1[1] - e2[0] > d):
                d =  e1[1] - e2[0] 
            if (e2[0] < e1[0] and e1[0] <= e2[1] and e2[1] < e1[1] and
                e2[1] - e1[0] > d):
                d = e2[1] - e1[0]
    return d

def setp_to_edges(p):
    """Transform a set partition into a list of edges.

    INPUT:
        p -- a Sage set partition.

    OUTPUT:
        A list of non-loop edges of the set partition. As this code just
        works with crossings, we can ignore the loops.

    EXAMPLE:
    The main example from the paper:
        sage: setp_to_edges(Set(map(Set, [[1,5],[2,4,9],[3],[6,12],[7,10,11],[8]])))
        [[7, 10], [10, 11], [2, 4], [4, 9], [1, 5], [6, 12]]
    """
    q = [ sorted(list(b)) for b in p ]
    ans = []
    for b in q:
        for n in range(len(b) - 1):
            ans.append(b[n:n+2])
    return ans
    
def dcrossvec_setp(n):
    """Return a list with the distribution of k-dcrossings on set partitions of [1..n].

    INPUT:
        n -- a nonnegative integer.

    OUTPUT:
        A list whose k'th entry is the number of set partitions p for
        which dcrossing(p) = k. For example, let L = dcrossvec_setp(3).
        We have L = [1, 0, 4]. L[0] is 1 because there's 1 partition of
        [1..3] that has 0-dcrossing: [(1, 2, 3)].

        One tricky bit is that noncrossing matchings get put at the end,
        because L[-1] is the last element of the list. Above, we have
        L[-1] = 4 because the other four set partitions are all
        d-noncrossing. Because of this, you should not think of the last
        element of the list as having index n-1, but rather -1.

    EXAMPLES:
        sage: dcrossvec_setp(3)
        [1, 0, 4]

        sage: dcrossvec_setp(4)
        [5, 1, 0, 9]

    The one set partition of 1 element is noncrossing, so the last
    element of the list is 1:
        sage: dcrossvec_setp(1)
        [1]
    """
    vec = [0] * n
    for p in SetPartitions(n):
        vec[dcrossing(setp_to_edges(p))] += 1
    return vec

def dcrossvec_cm(n):
    """Return a list with the distribution of k-dcrossings on complete matchings on n vertices.

    INPUT:
        n -- a nonnegative integer.

    OUTPUT:
        A list whose k'th entry is the number of complete matchings m
        for which dcrossing(m) = k. For example, let L =
        dcrossvec_cm(4). We have L = [0, 1, 0, 2]. L[1] is 1 because
        there's one matching on 4 vertices that is 1-dcrossing: [(2, 4),
        (1, 3)]. L[0] is zero because dcrossing() returns the *largest*
        k for which the matching has a dcrossing, and 0-dcrossing is
        equivalent to 1-dcrossing for complete matchings.

        One tricky bit is that noncrossing matchings get put at the end,
        because L[-1] is the last element of the list. Because of this, you
        should not think of the last element of the list as having index
        n-1, but rather -1.

        If n is negative, you get silly results. Don't use them in your
        next paper. :)

    EXAMPLES:
    The single complete matching on 2 vertices has no crossings, so the
    only nonzero entry of the list (the last entry) is 1:
        sage: dcrossvec_cm(2)
        [0, 1]

    Similarly, the empty matching has no crossings:
        sage: dcrossvec_cm(0)
        [1]

    For odd n, there are no complete matchings, so the list has all
    zeros:
        sage: dcrossvec_cm(5)
        [0, 0, 0, 0, 0]

        sage: dcrossvec_cm(4)
        [0, 1, 0, 2]
    """
    vec = [0] * max(n, 1)
    for m in CompleteMatchings(n):
        vec[dcrossing(m)] += 1
    return vec

def tablecolumn(n, k):
    """Return column n of Table 1 or 2 from the paper arxiv:0812.2725.

    INPUT:
        n -- positive integer.

        k -- integer for which table you want: Table 1 is complete
             matchings, Table 2 is set partitions.

    OUTPUT:
        The n'th column of the table as a list. This is basically just the
        partial sums of dcrossvec_{cm,setp}(n).

        table2column(1, 2) incorrectly returns [], instead of [1], but you
        probably don't need this function to work through n = 1.

    EXAMPLES:
    Complete matchings:
        sage: tablecolumn(2, 1)
        [1]

        sage: tablecolumn(6, 1)
        [5, 5, 11, 14, 15]

    Set partitions:
        sage: tablecolumn(5, 2)
        [21, 42, 51, 52]

        sage: tablecolumn(2, 2)
        [2]
    """
    if k == 1:
        v = dcrossvec_cm(n)
    else:
        v = dcrossvec_setp(n)
    i = v[-1]
    return [i + sum(v[:k]) for k in range(len(v) - 1)]
\end{lstlisting}

\bibliographystyle{amsplainurl} \bibliography{ref}
 
\end{document}